\documentclass[12pt]{article}
\usepackage{amsfonts,amsmath}
\newtheorem{theorem}{Theorem}

\begin{document}

\def\R{{\mathbb R}}
\def\C{{\mathbb C}}
\def\D{{\cal D}}
\def\u{{\bf u}}

\title {A fast decaying solution to the modified Novikov--Veselov equation with a one-point singularity}
\author{Iskander A. TAIMANOV
\thanks{Sobolev Institute of Mathematics, Academician Koptyug avenue 4, 630090, Novosibirsk, Russia, and Department of Mathematics and Mechanics, Novosibirsk State University, Pirogov street 2, 630090 Novosibirsk, Russia; e-mail: taimanov@math.nsc.ru. \newline
The work was supported by RSF (grant 14-11-00441). }}
\date{}
\maketitle

The modified Novikov--Veselov equation is a 2+1-dimensional soliton equation which was introduced in \cite{Bogdanov} and has the form
\begin{equation}
\label{mnv}
U_t = \big(U_{zzz} + 3U_z V + \frac{3}{2}UV_z \big) + \big(U_{\bar{z}\bar{z}\bar{z}} + 3U_{\bar{z}}\bar{V} + \frac{3}{2} U\bar{V}_{\bar{z}}\big),
\end{equation}
where
$$
V_{\bar{z}} = (U^2)_z,
$$
$z = x+iy \in \C$, $U$ is a real-valued function.
The function $V$ is uniquely defined under additional assumption: for instance, when considering fast decaying functions $U$ it is assumed that $V$ is also fast decaying.

This equation is represented by a Manakov $L,A,B$--triple:
$$
\D_t + [\D,A] - B\D = 0,
$$
where $\D$ is a two-dimensional Dirac operator:
$$
\D = \left(
\begin{array}{cc}
0 & \partial \\
-\bar{\partial} & 0
\end{array}
\right) + \left(
\begin{array}{cc}
U & 0 \\
0 & U
\end{array}
\right),
$$
with $\partial = \frac{\partial}{\partial z}$
and
$\bar{\partial} =
\frac{\partial}{\partial \bar{z}}$.
That implies that the equation preserves the ``zero energy level'', i.e. solutions to the equation
$\D\psi = 0$, deforming them via the equation $\frac{\partial \psi}{\partial t}= A\psi$.

If $U$ depends on $x$ only and $V = U^2$, then the modified Novikov--Veselov equation reduces to the modified Korteweg-de Vries equation
$$
U_ t = \frac{1}{4} U_{xxx}  + 6 U_x U^2.
$$
This fact it owes its name because the Novikov--Veselov equation \cite{NV}
is an analogous two-dimensional generalization of the Korteweg--de Vries equation.

The Cauchy problem, and its regularity for this equation were not studied until recently. In this paper we show that there are smooth initial data for the Cauchy problem such that
they decay as $O(r^{-2})$ as $r = \sqrt{x^2+y^2} = |z| \to \infty$ and the corresponding solution loses its regularity
 in finite time. Moreover the singularity of the solution which is extended for all times is supported at one point.

\begin{theorem}
The function
$$
U(x,y,t) =
-\frac{3((x^2+y^2+3)(x^2-y^2)-6x(C-t))}{Q(x,y,t)},
$$
\begin{equation}
\label{u1}
Q(x,y,t) = (x^2+y^2)^3 + 3(x^4+y^4)+18 x^2 y^2 +9(x^2+y^2) +
\end{equation}
$$
+ 9(C-t)^2 +
(6x^3-18xy^2-18x)(C-t),
$$
with $C = \mathrm{const}$,
satisfies the modified Novikov--Veselov equation.

Therewith

\begin{enumerate}
\item
the functions $U$ and $V$ are really-analytical everywhere outside the point
$x=y=0, t=C = \mathrm{const}$ at which they are not defined. The function $U$ has different finite limit values along the rays
$x/y = \mathrm{const}, t=C$, going into this point:
$$
\lim_{r \to 0, \varphi= \mathrm{const}} U(z,\bar{z},C) = - \cos\,2\varphi \ \ \ \ \mbox{for $z = r e^{i\varphi}$};
$$

\item
the functions $U$ and $V$ decay as $O(r^{-2})$ as $r \to \infty$ for all fixed times;

\item
the first integral (conservation law) $\int_{\R^2} U^2 dx\, dy$ of the equation (\ref{mnv})
is equal to  $3\pi$ for $t \neq C$ and jumps to $2\pi$ at $t=C$.
\end{enumerate}
\end{theorem}

Let us expose for completeness the formula for $V(x,y,t)$:
$$
V =
\frac{(z(\bar{\gamma}-\gamma)-\delta z^2 -\bar{\delta})^2}{(|\gamma|^2+|\delta|^2)^2} +
\frac{2U}{1+|z|^2}  - 2 \frac{iz(z(\bar{\gamma}-\gamma)-\delta z^2 -
\bar{\delta})}{(|\gamma|^2+|\delta|^2)(1+|z|^2)},
$$
where
$$
\gamma = i(x^2-y^2),
$$
$$
\delta = y\left(1+x^2 - \frac{y^2}{3}\right) -i\left[x\left(1+y^2-\frac{x^2}{3}\right)+C-t\right].
$$

This solution is constructed by using the Moutard transformation for two-dimensional Dirac operators \cite{C}. Given a solution $\psi_0$ to the equation $\D\psi_0 = 0$, the
transformation gives a new Dirac operator $\D_M$ with another potential and supplies an analytical procedure for finding all solutions of the equation
$\D_M\varphi=0$ from solutions to the equation $\D\psi = 0$. This transformation is expanded to a transformation of
solutions to the modified Novikov--Veselov equation.

The explicit formulas for the Moutard transformation are quite bulky and we skip them. In \cite{T1} it is shown that this transformation has a geometrical interpretation
in terms of conformal (M\"obius) geometry and the Weierstrass representation of surfaces \cite{T-RS}.

It is with the help of this interpretation the solution  (\ref{u1}) was obtained.
Let us explain its geometrical meaning.

The Enneper surface is a minimal surface which is immersed into $\R^3$. Up to translations it is defined by the formulas
$$
u^1(x,y) = y \left(\frac{y^2}{3}-x^2 -1\right) + u^1_0,
$$
$$
u^2(x,y) = x \left(1+ y^2 - \frac{x^2}{3}\right) + u^2_0,
$$
$$
u^3(x,y) = x^2-y^2 + u^3_0,
$$
where $u^1,u^2,u^3$ are the Euclidean coordinates in $\R^3$ and
${\bf u}_0 = (u^1_0,u^2_0,u^3_0)$ is the image of the origin $x=y=0$
under the immersion. The induced metric (the first fundamental form) is conformally Euclidean:
$ds^2 = g_0(x,y) dz\,d\bar{z}$.

Let us consider the family $S_t$ of translated Enneper surfaces corresponding to ${\bf u}_0 = (0,C-t,0)$. The inversion of $\R^3$ centered at the origin has the form
$$
T: \u \to -\frac{\u}{|\u|^2}, \ \ \ \u = (u^1,u^2,u^3) \in \R^3.
$$
Let us correspond to every surface $S_t$ its image $\Sigma_t$ under the inversion. Since the inversion is a conformal mapping, the
induced metrics on $\Sigma_t$ are also conformally Euclidean:
$ds^2 = g(x,y,t) dz\,d\bar{z}$.
Let us define the function
$$
U(x,y,t) = \frac{H \sqrt{g}}{2},
$$
where $H$ is the mean curvature of $\Sigma_t$ at $(x,y)$ and
$g$ is the conformal factor of the metric at the same point. This is the function (\ref{u1}).
It has singularities only at points at which $S_t$ passes through the origin and therefore which are mapped into the infinity by the inversion.
It is easy to check that there is only one such a point: $x=y=0, t=C$.

The detailed exposition of these results will be given elsewhere.

Let us remark that by using the Moutard transformation of two-dimen\-sio\-nal Schr\"odinger operators blowing up solutions of the
Novikov--Veselov equation were constructed in
в \cite{TT1,TT2}. However they do not have such a geometrical interpretation.

\end{document}